\documentclass[12pt]{amsart}

\usepackage{amssymb}
\usepackage{amsmath}
\usepackage{amsthm}
\usepackage{amssymb}
\usepackage{marginnote}

\theoremstyle{plain}
\newtheorem{theorem}[subsubsection]{Theorem}

\newtheorem*{theorem*}{Theorem}
\newtheorem{proposition}[subsubsection]{Proposition}
\newtheorem*{proposition*}{Proposition}

\newtheorem*{lemma*}{Lemma}
\newtheorem{corollary}[subsubsection]{Corollary}
\newtheorem*{corollary*}{Corollary}

\theoremstyle{definition}

\theoremstyle{remark}

\newtheorem*{eexample*}{Example}

\newcommand{\0}{{(0)}}

\newcommand{\CC}{\mathbb{C}}

\newcommand{\RR}{\mathbb{R}}

\newcommand{\ZZ}{\mathbb{Z}}

\newcommand{\OO}{\mathcal{O}}
\newcommand{\WW}{\mathcal{W}}

\newcommand{\Hom}{\operatorname{Hom}}

\newcommand{\meas}{\operatorname{\meas}}

\newcommand{\nc}{\newcommand}
\nc{\G}{\Gamma}
 \nc{\sm}{\setminus}
 \nc{\sub}{\subset}
 \nc{\lm}{\lambda}
  \nc{\Lm}{\Lambda}
 \nc{\al}{\alpha}
 \nc{\bt}{\beta}
 \nc{\om}{\omega}
 \nc{\dl}{\delta}
 \nc{\g}{\gamma}
 \nc{\Dl}{\Delta}
 \nc{\Om}{\Omega}
 \nc{\s}{\sigma}
 \nc{\ro}{\rho}
 \nc{\te}{\theta}
 \nc{\SLR}{\operatorname{SL}_2(\RR)}
 \nc{\GLR}{\operatorname{GL}_2(\RR)}
 \nc{\PGLR}{\operatorname{PGL}_2(\RR)}
 \nc{\PSLR}{\operatorname{PSL}_2(\RR)}
 \nc{\PSLZ}{\operatorname{PSL}_2(\ZZ)}
 \nc{\SLC}{\operatorname{SL}_2(\CC)}
 \nc{\uH}{\mathbb H}
 \nc{\fD}{\mathcal{D}}
 \nc{\fE}{\mathcal{E}}
 \nc{\fO}{\mathcal{O}}
 \nc{\haf}{\frac{1}{2}}
 \nc{\qtr}{\frac{1}{4}}
 \nc{\shaf}{{\scriptstyle\frac{1}{2}}}
 \nc{\hlm}{{\scriptstyle\frac{\lambda}{2}}}

 \nc{\8}{\infty}
 \nc{\7}{{-\infty}}
 \nc{\inv}{^{-1}}
 \nc{\eps}{\varepsilon}
 \nc{\aG}{\mathbf{G}}
 \nc{\spn}{\operatorname{Span}}
 \nc{\Cm}{\operatorname{CM}}

\nc{\Ck}{\textsl{k}}

\nc{\fa}{\mathfrak a} \nc{\fg}{\mathfrak g} \nc{\fii}{\mathfrak i}\nc{\fk}{\mathfrak k}
\nc{\fh}{\mathfrak h} \nc{\fm}{\mathfrak m} \nc{\fn}{\mathfrak n}
\nc{\fA}{\mathfrak A} \nc{\fC}{\mathfrak C} \nc{\fI}{\mathfrak I}
\nc{\fL}{\mathfrak L} \nc{\fS}{\mathfrak S}
\nc{\fz}{\mathfrak z} \nc{\fl}{\mathfrak l}
\nc{\fp}{\mathfrak p}
\nc{\ft}{\mathfrak t}
\makeatletter
\@addtoreset{equation}{section}
\makeatother
\numberwithin{equation}{section}

\begin{document}
\title[Lattice points and period bounds]{\bf Lattice points counting and bounds on periods of Maass forms}

\author{Andre Reznikov} \email{reznikov at math.biu.ac.il}


\author{Feng Su} \email{fsu at math.biu.ac.il}

\address{Department of Mathematics\\
Bar-Ilan University \\
Ramat Gan, 52900 \\
Israel}
\thanks{  The research was partially supported  by the ERC grant  291612 and by the ISF grant 533/14.}

\allowdisplaybreaks
\vspace*{-0.8cm}
\begin{abstract} We provide a ``soft'' proof for non-trivial bounds on spherical, hyperbolic and unipotent Fourier coefficients of a fixed Maass form for a general co-finite  lattice $\G$ in $\PGLR$. We use the amplification method based on the Airy type phenomenon for corresponding matrix coefficients and an effective Selberg type pointwise asymptotic for the lattice points counting in various homogeneous spaces for $\PGLR$. This requires only $L^2$ theory.  We also show how to use the uniform bound for the $L^4$-norm of $K$-types in a fixed automorphic representation of $\PGLR$ obtained in \cite{br} in order to slightly improve these bounds. 
\end{abstract}

\maketitle

\section{Introduction}\label{intr}
In this  paper,  we explore an interplay between two classical topics in analytic theory of automorphic functions. Namely, we show how to use lattice points counting in various $\PGLR$-homogeneous spaces   in order to obtain non-trivial bounds on Fourier coefficients of Maass forms (or more generally, on generalized periods of Maass forms). We start with an  analysis of certain (generalized) matrix coefficients of unitary irreducible representations of $\PGLR$. It is a well-known fact that such matrix coefficients  have well-understood and uniform asymptotic behavior at infinity.  Our point of the departure is that, in spite of this   fact,  in certain ranges of parameters (i.e., parameters of the representation and the group variable) matrix coefficients might exhibit ``abnormally'' big values. The origin of this phenomenon is  classical  and  best associated with the Airy function (see \cite{ar}). We use this property of matrix coefficients in order to ``amplify'' (in the language of analytic number theory) certain automorphic coefficients. A somewhat unexpected fact is that the Airy type phenomenon for  matrix coefficients holds (for appropriate values of parameters)  over domains in $\PGLR$ with a {\it large} volume.  This allows us to use as a  global input  the {\it effective} lattice point counting  in corresponding  domains in $\PGLR$ or in relevant homogeneous spaces (i.e., the asymptotic  count of lattice points of a discrete subgroup  with a polynomial saving in the remainder).

\subsection{Maass forms}\label{M-forms}
Let $Y$ be a compact Riemann surface  with a Riemannian metric of
constant curvature $-1$. We denote by $dv$ the associated volume
element and by ${\rm d}(\cdot,\cdot)$ the corresponding distance
function. The corresponding Laplace-Beltrami operator $\Dl$ is
non-negative and has purely discrete spectrum on the space
$L^2(Y,dv)$ of functions on $Y$. We will denote by $0=\mu_0< \mu_1
\leq \mu_2 \leq ...$ the eigenvalues of $\Dl$ and by $\phi_i$ the
corresponding eigenfunctions normalized by the $L^2$-norm on $Y$.
In the theory of automorphic functions, the functions $\phi_i$ are
called non-holomorphic forms, Maass forms (after Hans Maass) or simply automorphic functions on $Y$. We will write the eigenvalue of a Maass form   $\phi_i$ in the form  $\mu_i=(1-\lm_i^2)/4$ with $\lm_i\in i\RR\cup (0,1]$. 

The study of Maass forms is important in  number theory, analysis
and in  mathematical physics. In
particular, various questions concerning analytic properties of
 eigenfunctions $\phi_i$ have drawn much  attention in recent
years (see surveys \cite{Sa2}, \cite{Sa3}).
In this paper, we study analytic properties of periods of
eigenfunctions over some special curves on $Y$ (see  \cite{Z} for a survey of recent analytic developments). 

A function on $Y$ could be lifted to the universal cover which is the hyperbolic upper half-plane $\uH=G/K$ for $G=\PGLR$ and $K=PO(2)$. We have $Y=\G\sm\uH$ with $\G=\pi_1(Y)\subset \PGLR$.  When convenient, we will view Maass forms $\phi_i$ as $\G$-invariant eigenfunctions on $\uH$ or on $G$. 

For completeness, we describe now in some detail spherical Fourier expansion  on $Y$ in classical terms. This expansion could be viewed as a substitute for the Taylor expansion at a point in $Y$. Two other types of Fourier expansions, the unipotent and the hyperbolic ones, are discussed at length in  \cite{re1} and we refer to the discussion there of these expansions in classical terms. We present all three expansions in representation theoretic terms later in the paper. 

\subsubsection{Spherical Fourier coefficients} Let $\phi_i$ be a Maass form with the eigenvalue $\mu_i>0$. 
Writing the Laplace-Beltrami operator in polar
geodesic coordinates $(r,\te)$, $r\geq 0$, $1> \te\geq 0$,  centered at a point $y\in Y$, and using the separation of
variables, one sees immediately that there exists a special function
$C_\mu(r,n)$ such that  $\phi_i$ (or, what is the same, the  lift of $\phi_i$ to $\uH$)
has the Fourier series expansion given by
\begin{align}\label{c-fourier}
\phi_i((r,\te))=\sum_n a_n(y,\phi_i)C_{\mu_i}(r,n)e^{2\pi
in\te}\ , \ n\in \ZZ .
\end{align}
This is one of the classical expansions considered in the theory of automorphic functions (see \cite{pe}, \cite{Good}, \cite{he}, \cite{Sa1}). 
The function $C_\mu(r,n)$ is unique up to a multiplicative constant (which classically is  chosen by fixing an explicit integral representation or an asymptotic expansion for $C_\mu(r,n)$), and  depends only on the
eigenvalue $\mu_i$ of $\phi_i$ and not on the choice of the eigenfunction $\phi_i$. The function $C_\mu(r,n)$
 is essentially equal to the appropriate hypergeometric
function or the Legendre function. We will relate  $C_\mu(r,n)$ to a matrix coefficient from representation theory (this is gives one of the classical integral representations of $C_\mu(r,n)$). 
The expansion \eqref{c-fourier} is valid in any disc around $y$ (and, in fact, on the whole upper-half plane uniformizing $Y$). The expansion \eqref{c-fourier} is similar
to the Taylor expansion at the point $y\in Y$, and in fact would be the Taylor expansion for the holomorphic forms (see \cite{Good}).

It is important for us that the coefficient $a_n(y,\phi_i)$
is given by a (generalized) period of a Maass form. Namely, denoting a geodesic circle centered at a point $y\in Y$ by  $\s(r,y)=\{y'\in Y|{\rm
d}(y',y)=r\}$. We have 
\begin{align}\label{F-coef-K}
a_n(\phi_i,y)C_{\mu_i}(r,n)=\int\limits_{\te\in\s(r,y)}\phi_i(\te)e^{-2\pi
in\te}d\te .
\end{align} We stress that the function $C_\mu(r,n)$ is independent of the point $y$ and depends only on the eigenvalue, while
coefficients $a_n(\phi_i,y)$ capture the structure of the
eigenfunction $\phi_i$.

The coefficients $a_n(\phi_i,y)$ are important objects in Number Theory and were studied extensively in recent years owing to their connection to $L$-functions via the celebrated theorem of J.-L. Waldspurger \cite{wa} (for an arithmetic $Y$, special points $y$ and for so-called Hecke-Maass forms).

We are interested in bounds on these coefficients. In \cite{re1}, we obtained the following mean-value bound (in fact, a similar bound for unipotent Fourier coefficients  is the classical bound of E. Hecke/G. Hardy).  For any $y\in Y$, there
exists a constant $C_{y}>0$ such that for any $T\geq 1$, the
following bound holds:
\begin{align}\label{a-bound-1}
\sum_{ |n|\leq T} |a_n(\phi_i,y)|^2\leq C_{y}\cdot \max(T,\mu_i^\haf)\ ,
\end{align} for all $\phi_i$. The bound is sharp in the range $T>\mu_i^\haf$ and its improvement in the range $\mu^\haf\ll T$ is one of the major problems in the field (it is expected that $|a_n(\phi_i,y)|\ll \mu_i^\eps$). 

In this paper we describe a simple argument which allows one  to improve such a bound for a {\it fixed} Maass form. Namely, we will prove the following bound on the average size of coefficients $a_n$  in ``short intervals'' (i.e., in intervals of the type $[T,T+T^\theta]$ with $\theta<1$)
\begin{equation}\label{sph-thm-bound-intro}
\sum\limits_{|n-T|\leq T^{5/7}}|a_n(\phi_i,y)|^2\leq C'_{y,\phi_i}\cdot T^{6/7}\ .
\end{equation} Note that we miss the sharp on average bound in the above interval, the so-called Lindel\"{o}f on average bound. Within $L^2$ theory alone, we obtain the exponent $8/9$ instead of $6/7$ above. We also prove similar {\it non-trivial} bounds on hyperbolic and unipotent Fourier coefficients of Maass forms by the same technique. 

 Bounds for coefficients $a_n(\phi_i,y)$ were investigated  in recent years with non-trivial (i.e., better than $n^\haf$) bounds obtained in \cite{v}, \cite{re2} and \cite{bl-etc}. So far, the best  bound  for a general point $y$ and/or general $\phi_i$, was obtained in \cite{re2} and reads  $\sum\limits_{|T-n|\leq T^{2/3}}|c_n|^2\ll T^{2/3+\eps}$.   This is on par with H. Weyl/G.~Hardy-J. Littlewood bound for the Riemann zeta function and also  ``Lindel\"{o}f on average'' type of a bound in short intervals.  Any polynomial improvement over the trivial bound $|n|^\haf$ leads to a subconvexity bound for the corresponding $L$-function. In analogy with the Ramanujan/Lindel\"{o}f conjecture one expects the bound $|a_n|\ll n^\eps$ to hold for any $\eps>0$. 

Our point in the present paper (which owns main ideas to \cite{br2} and to \cite{bl-etc}) is not in the improvement of the exponent in the bound \eqref{sph-thm-bound-intro} (which we fail to achieve!), but in the simplicity of the method and its relation to the problem of counting  lattice points. In particular, we do not use Rankin-Selberg type integrals (or, which is essentially the same, triple products from  \cite{re2}). 

{\bf Acknowledgments.} This paper is a part of a joint project
with Joseph Bernstein whom we would like to thank for numerous fruitful
discussions. It is also a pleasure to thank Amos Nevo for explanations concerning lattice points counting.

The research was partially supported  by the ERC grant  291612 and by the ISF grant 533/14.

\section{Representation theory}\label{rep}  
In this section, we recall some basic representation theory of $G=\PGLR$ (see the classic \cite{ge} and \cite{br} for  notations we will use).
We will denote elements in the unipotent subgroup $N$ by $n_x=\left(\begin{smallmatrix}1&x\\0&1
\end{smallmatrix}\right)$, $x\in\RR$, and in  the maximal compact subgroup $K=PO(2,\RR)$ by  $k_{\theta}=\left(\begin{smallmatrix}
\cos\theta&\sin\theta\\-\sin\theta
&\cos\theta\end{smallmatrix}\right)$, $\theta\in\RR/2\pi\ZZ$. 
We will often construct elements in $G=\PGLR$ by choosing their representatives in $\GLR$.

 We will work with smooth irreducible representations. Moreover, we will restrict ourselves to spherical representations (i.e., those with a non-zero $K$-fixed vector) with the trivial central character. There are even and odd spherical representations of $G$ with respect to the action of elements with negative determinant (e.g., with respect to the action of $\dl=\left(\begin{smallmatrix}
-1  & \\
&\ \ 1 \\ \end{smallmatrix}\right)$ on a $K$-fixed vector; see Sect. 2.5, \cite{bu}, \cite{re1}).  We will deal with even representations and the treatment of odd ones is identical. Spherical irreducible representations have various geometric realizations, called models. The basic model is the model in the space of homogeneous functions on the plane.  
For $\lambda\in\mathbb{C}$, we
denote by $W_{\lambda}$
the space of smooth even homogeneous functions on $\mathbb{R}^2\smallsetminus\{(0,0)\}$
of the homogeneous degree $\lambda-1$, i.e., $f(ax,ay)=|a|^{\lambda-1}f(x,y)$ for any $a\in\mathbb{R}\smallsetminus\{0\}$.
An element $g\in \GLR$ acts on $f\in W_{\lambda}$
via $\pi_{\lambda}(g)f(x,y)=
|\det g|^{(\lambda-1)/2}f\left(g^{-1}(x,y)\right)$ where $g^{-1}(x,y)$ means the usual matrix multiplication on columns. 
This is a well-defined action of $G$ in  view of the homogeneity of $f$ (i.e., it is trivial on the center of $\GLR$).
The nontrivial unitary spherical dual of $G$ is given
by $(\pi_{\lambda},W_{\lambda})$ for 
$\lambda\in(-1,1)\cup i\,\mathbb{R}$. A homogeneous function $f\in W_\lm$ is uniquely determined by its values  
on a line in $\RR^2\sm 0$ (e.g., on the line $L_x=\{(x,1)\ |\ x\in \RR\}\subset \RR^2$).  We denote the space of the corresponding model by $V_{\lambda}$ and called it the {\it line model}.
We still   use the notation $\pi_{\lambda}$ to denote the action of $G$ on $V_{\lambda}$. We will assume for simplicity that $\lm\in i\RR$, i.e., $\pi_\lm $ is a representation of the {\it principal series}. In terms of the eigenvalue of a Maass form, this means that $\mu\geq 1/4$. We denote by  $\langle \cdot,\cdot\rangle=\langle \cdot,\cdot\rangle_{V_\lm}$ the invariant Hermitian form in the representation $V_\lm$. In the line model, it is given by 
$\langle v,u\rangle=\pi\inv \int_\RR v(t)\bar u(t)dt$. 

We denote by $X=\G\sm G$ the automorphic space and by 
$dx$  the $G$-invariant Radon measure on $\Gamma\backslash G$ induced by the Haar measure $dg$ of $G$. Co-finiteness of $\G$ means that the total volume of $X$ is finite and we can normalize $dx$ to have the total volume one. We will assume for simplicity that  $\dl=\left(\begin{smallmatrix}
-1  & \\
&\ \ 1 \\ \end{smallmatrix}\right)\in \G$, and in particular the space $X=\G\sm G$ has one connected component.

The (right) action of $G$ on $L^2(X)$ is given by $R(g)\phi(x)=\phi(xg)$, $x\in X$. An automorphic representation is a pair $(\nu,\pi)$, where $\pi$ is an abstract irreducible unitary representation of $G$ in a space $V_\pi$ of smooth vectors and $\nu:V_\pi\to L^2(X)$ is a $G$-homomorphism which we assume to be an isometry. For a vector $v\in V_{\pi}$, we denote by $\phi_v$ the corresponding automorphic function $\phi_v=\nu(v)\in  L^2(X)$. It is well-known that for a smooth vector $v\in V_\pi$, $\phi_v$ is a smooth function on $X$. In particular, restriction of such automorphic functions to subsets in $X$ is well-defined.

\section{Spherical periods}
In this section we consider  spherical periods, i.e., those coming from orbits of the compact subgroup $K$.
\subsection{Equivariant functionals}\label{form}
We fix an automorphic representation $(\pi,\nu)$. We will assume that $(\pi, V_\pi)$ is an even representation of the principal series  for some $\lm\in i\RR$.  Let $v\in V_\lm$ be a smooth vector. We consider the corresponding (smooth) automorphic function $\phi_v$. Let ${\OO}=x_0K\subset X$, $x_0=\G g_0\in X$, be a $K$-orbit, $\OO\simeq Stab_K(x_0)\sm K$. The group  $Stab_K(x_0)\subset K$ is finite and is non-trivial only if $x_0$ is a fixed point on $X$ of an elliptic element in $\G$. We denote by $i_\OO=|Stab_K(x_0)|$ the order of this group.   We will fix a $K$-invariant measure $do$ on $\OO$ and will assume for simplicity that $vol_{do}(\OO)=1$. We can restrict $\phi_v$ to $\OO$ and expand it as $$\phi_v|_\OO=\sum_na_n(v)\chi^._n,\quad a_n(v)\in\mathbb{C}$$ where $\chi_n^.$'s are functions on $\OO$ corresponding to unitary characters $\chi_n: K\to S^1$ which are trivial on $Stab_K(x_0)$ (i.e., $\chi^._n(x_0k)=\chi_n(k)$). Functions $\chi^._n$ constitute an orthonormal basis of $L^2(\OO)$. 

Coefficients $a_n(v)$ are given by an automorphic period $$a_n(v)=\int_{o\in\OO}\phi_v(o)
\overline{\chi^._n(o)}do\ .$$
Such an integral defines a $\chi_n$-equivariant functional $\ell^{\rm aut}_n\in\Hom_K(V_\pi,\CC)$:
$$\ell^{\rm aut}_n\,:\,V_\pi\rightarrow\mathbb{C},\quad v\mapsto\int_{\OO}\phi_v(o)
\overline{\chi^._n(o)}do.$$
It is known that the space of such functionals is at most one-dimensional. This is equivalent to the fact that all $K$-types are at most one-dimensional for the group $\PGLR$. We can consider the   $\chi_n$-equivariant functional on the abstract representation $V_\pi\simeq V_{\lm}$ which corresponds to a normalized $K$-type. We call such a functional the model functional. Namely, let $e_n\in V_\lm$ be a vector of $K$-type of weight $n$, $||e_n||_{V_\pi}=1$,  and 
$$\ell^{\rm mod}_n\,:\,V_{\lambda}\rightarrow\mathbb{C},\quad v\mapsto\langle v,e_{n}\rangle$$
the corresponding $\chi_n$-equivariant functional.  The functional $\ell^{\rm mod}_n$ is well-defined up to a constant of the absolute value one. A simple computation shows that a normalized $K$-type of weight $n$ in the line model $V_{\lambda}$ is given by $e_n(x)=e^{in\arctan x}(1+x^2)^{\frac{\lambda-1}{2}}$ and  the functional $\ell^{\rm mod}_n$ is given by the integral $$\ell^{\rm mod}_n(v)={\pi}\inv\int_{\mathbb{R}}v(x)\overline{e_n(x)}dx\ .$$ 
In view of the multiplicity one $\dim\Hom_K(V_\pi,\chi_n)\leq 1$ of $K$-equivariant functionals,  there exists a scalar $c_n=c_n(\OO,\nu)\in\mathbb{C}$ such that 
\begin{equation}\label{c-n}
\ell^{\rm aut}_n(\phi_v)=
c_n\cdot\ell^{\rm mod}_n(v)\ ,
\end{equation}
for all $v\in V_{\lambda}$.
Note that coefficients $c_n$ depend only on $\OO$, $\nu$ and $n$, but not on 
 a vector $v$. 

From geometric considerations, one obtains a very general bound for these coefficients (sometimes called a convexity bound or Hecke/Hardy bound). Namely, one has the following  bound on the average size of coefficients $c_n$ over ``long intervals'' $$\sum\limits_{|n|\leq T}|c_n|^2\leq C_\G \max(|\lm|, T)\ .$$ This bound is sharp for $T\geq |\lm|$ (see \cite{re1}).   

Our main result in this section is the following
\begin{theorem}\label{sph-thm} Let $\G\subset G$ be a co-compact lattice, $\OO\subset X$ be a $K$-orbit and $(\pi,\nu)$ be an automorphic representation of the principal series. We have the following bound for the spherical Fourier coefficients $c_n\in\CC$  
\begin{equation}\label{sph-thm-bound}
|c_n|^2\leq C (1+|n|)^{6/7}\ ,
\end{equation}
for some constant $C>0$ depending on $\G$, $\OO$ and $(\pi,\nu)$, but not on $n$. 
\end{theorem}
In fact, we will prove the following bound on the average size of coefficients $c_n$  in ``short intervals''
\begin{equation}\label{sph-thm-bound-aver}
\sum\limits_{|T-n|\leq T^{5/7}}|c_n|^2\leq C\cdot T^{6/7}\ .
\end{equation}

Staying only within the $L^2$ theory (i.e., without the use of the uniform bounded on $L^4$-norm of $K$-types; see Theorem \ref{L4-norm}), we obtain a slightly weaker bound $\sum\limits_{|T-n|\leq T^{7/9}}|c_n|^2\leq C'\cdot T^{8/9}$.  

We also discuss non-uniform lattices in Section \ref{non-comp}.

As we mentioned in the Introduction, better bounds on coefficients $c_n$ are known  by more sophisticated methods.

\subsubsection {Integrated Hermitian form} The proof is based on the by now  standard idea of ``fattening the cycle'', i.e., averaging over the group action the Hermitian form coming from the $L^2$-form on the orbit $\OO$. This technique was introduced (among other places) in \cite{br2} and in \cite{v}. The new ingredient in the present work is the use of the Airy amplification of matrix coefficients (in a slightly different context this idea appeared before  in  \cite{br}) together with averaging over {\it big thin shells} (this is inspired by \cite{bl-etc} where effective equidistribution over a shifted cycle is used).    

By the Plancherel formula, the $L^2$-norm of the restriction $\phi_v|_\OO$ can be written  as 
$$
\|\phi_v\|^2_{L^2(\OO)}=\sum\limits_n|a_n(v)|^2=\sum\limits_n\left|\ell^{\rm aut}_n(\phi_v)\right|^2
=\sum\limits_n|c_n|^2\cdot
\left|\ell^{\rm mod}_n(v)\right|^2.$$
We can act on this identity by the action of $G$ on the vector $v$. For  $g=ak\in AK$ (in fact we can act by a general $g\in G$, but the action by $K$ on the ``left'' is clearly trivial), we have
\begin{equation}\label{va}
\|\phi_{\pi_{\lambda}(ak)v}\|^2_{L^2(\OO)}=\sum\limits_n|c_n|^2\cdot
\left|\ell^{\rm mod}_n\big(\pi_{\lambda}(ak)v\big)\right|^2.
\end{equation}
We now want to integrate this identity with respect to the action of $G$. Let $\al\in C^\8_0(G//K)$ be a smooth real valued non-negative  function on $G$
such that 
\begin{itemize}
\item
$\al(g)\leq 1$;
\item
$\al$ has a compact support;
\item 
$\al$ is bi-$K$-invariant.
\end{itemize}
The exact form of the function $\al$ will be specified later (see Section \ref{test-al}). We use $KA^-K$ decomposition instead of the customary $KA^+K$ decomposition for  later convenience (here $A^\pm\subset A$ as usual denote sets of diagonal matrices $diag(t,t\inv)$ with $t\geq 1$ or $t\leq 1$ correspondingly). For $a\in\mathbb{R}_+$ (in what follows $a>1$), we denote by $[a]$ (or by $a$ when it does not lead to  confusion) the diagonal matrix ${\rm diag}(a^{-1},a)\in A$ ($a>1$). 
We multiply by $\al$  the  identity \eqref{va} and
integrate it over $ak\in A^-K$ with respect to the measure $\rho(a)d[a]dk=\tfrac{a^{-2}-a^2}{2}\frac{da}{a}dk$ coming from the Cartan decomposition (see \cite[Proposition 5.28]{kn}). We normalize the measure $dk$ on $K$  so that ${\rm vol}(K)=1$, and hence we have $\int_\OO f (o)do=i_\OO\inv\int_K f(x_0k)dk$ for any $f\in L^1(\OO,do)$ (here $i_\OO=|Stab_K(x_0)|$).  We obtain on the left in \eqref{va} 
\begin{eqnarray}\label{fold-sph} \nonumber
&&\int_{A^-K}\left(\int_\OO|\phi_{\pi_{\lambda}(ak)v}(o)|^2do\right)\al(ak)\rho(a)dadk \\
&&=\int_{A^-K}\left(i_\OO\inv\int_K|\phi_{\pi_{\lambda}(ak_2)v}(x_0k_1)|^2dk_1\right)\al(ak_2)\rho(a)dadk_2 \notag\\\nonumber
&&=i_\OO\inv\int_{A^-K}\left(\int_K|R(ak_2)\phi_v(x_0k_1)|^2dk_1\right)\al(ak_2)\rho(a)dadk_2\\\nonumber
&&=i_\OO\inv\int_{A^-K}\int_K|\phi_v(x_0k_1ak_2)|^2\al(ak_2)\rho(a)dk_1dadk_2 \\ \nonumber
&&\overset{(\ast)}{=}\notag i_\OO\inv\int_G|\phi_v(x_0g)|^2\al(g)dg=
 i_\OO\inv\int_G|\phi_v(g)|^2\al(g_0\inv g)dg \\ &&
=i_\OO\inv\int_{\Gamma\backslash G}|\phi_v(x)|^2 H_{\al,x_0}(x)dx \ ,
\end{eqnarray} where $H_{\al,x_0}(x)=\sum_{\gamma\in\Gamma}
\al(g_0\inv\gamma g)\in C(X)$ with $x_0=\G g_0\in X$ and $x=\G g\in X$.  In the step ($\ast$) we have used the integral formula for the Cartan $KAK$ decomposition.
On the right hand side of \eqref{va} we get the integral of the matrix coefficient. Hence we obtain the following integrated identity

\begin{equation}\label{eq}
i_\OO\inv\int_{\Gamma\backslash G}
|\phi_v(x)|^2 H_{\al,x_0}(x)dx=
\sum\limits_n|c_n|^2\cdot\int_{A^-K}
\left|\langle\pi_{\lambda}(ak)v,
e_n\rangle\right|^2\al(ak)\rho(a)dadk.
\end{equation}

We will apply this identity to a $K$-type vector $v=e_m$, $m\in\ZZ$. We have $\pi_{\lambda}(k_{\theta})e_m=e^{im\theta}e_m$, and hence the right hand side of (\ref{eq}) simplifies  to 
\begin{equation}\label{for-em}\sum\limits_n|c_n|^2\cdot\int_{A^-}
\left|\langle\pi_{\lambda}(a)e_m,
e_n\rangle\right|^2\al(a)\rho(a)da\ .
\end{equation}

\subsection{Airy asymptotic}\label{airy} We now study the asymptotic of the matrix coefficient $\langle\pi_{\lambda}([a])e_m,
e_n\rangle$ as a function of $m$, $n$ and $a$ for  $\lm\in i\RR$ which is {\it fixed}. We find a range of $m$, $n$ and $a$ where the value of the matrix coefficient is ``abnormally'' big due to the Airy type phenomenon, and use this to amplify a particular range of $n$'s in the identity \eqref{eq}. 

  It is more convenient to deal with the matrix coefficient $\langle e_m,\pi_{\lambda}
\left([a]\inv\right)
e_n\rangle =\langle\pi_{\lambda}
\left([a]\right)e_m,
e_n\rangle$. 
We have $\pi_{\lambda}\left([a]\inv\right)
e_n(x)=|a|^{\lambda-1}e_n(a^{-2}x)$ ($[a]\inv=diag(a, a\inv)$) and
\begin{eqnarray}\label{1}
&&\langle e_m,\pi_{\lambda}
\left([a]\inv\right)
e_n\rangle=\\ \nonumber &&\pi\inv |a|^{\lm-1}
\int_{\mathbb{R}}e^{im\arctan (x)-in
\arctan (x/a^2)}(1+x^2)^{\frac{\lambda-1}{2}}(1+(x/a^2)^2)^{\frac{-\lambda-1}{2}}dx\ .
\end{eqnarray}
We study the  integral in \eqref{1} by the method of the stationary phase method  viewing $n$ as a large parameter and $m$ and $a$ as auxiliary parameters. We  give an 
asymptotic expression for the  integral in \eqref{1} in terms of the classical Airy function (see \cite{mo}, \cite{ol}). To study the Airy stationary point near zero, we introduce new parameters  $\beta$ and $\varepsilon$ such that  $n=\beta^2m$ and $a=\beta(1+\varepsilon)$. We split (with  help of a standard smooth partition of unity by functions $\chi_0$ supported near $0$ and $\chi_\8$ supported off $0$, $\chi_0+\chi_\8\equiv 1$)  the integral into two integrals over $[-1, 1]$ and the complement. A simple computation shows that in the oscillating integral over the complement of $[-1,1]$ the phase function has no critical points (including at infinity) and hence the integral is small as $n\to\8$ (i.e., decays faster than $n^{-N}$ for any $N\geq 1$). We are left to deal with the integral over $[-1,1]$.

Denote by $S(\bt,\eps;x)=\arctan(x)-\bt^{2}\arctan (\bt^{-2}(1+\eps)^{-2}x)$ the phase of the oscillating integral \eqref{1}
and by $f(\bt,\eps;x)=(1+x^2)^{\frac{\lambda-1}{2}}(1+(\bt(1+\eps))^{-4}x^2)^{\frac{-\lambda-1}{2}}$ its amplitude.  The relevant part of the integral in (\ref{1}) could be written as  
\begin{equation}\label{I-mna} I_{m,n}(a):=\int_{\mathbb{R}}
e^{imS(\bt,\eps;x)}f(\bt,\eps;x)\chi_0(x)dx\ .
\end{equation}
 The phase function has an Airy type critical point at $x=0$. We have the classical asymptotic $\int_{\RR}
e^{it(\dl y+y^3/3)}f(y)dy\sim C t^{-1/3}{\rm Ai}(t^{2/3}\dl)+O(t^{-2/3})$ with some explicit constant $C>0$ and for $f\equiv 1$ near $0$, smooth and of a compact support (see \cite{mo}). The Airy function satisfies $|{\rm Ai}(t)|\geq c>0$ for $|t|\leq 1$ (e.g., the first zero of ${\rm Ai}$ satisfies $|t_1|\approx 2.388$). We have the following expansion: 
\begin{equation}\label{S} S(\bt,\eps;x)=2\eps x-x^3/3+({\it smaller \ order\ terms \ in \ } x, \eps, \bt\inv)\ .
\end{equation}
This suggests  that the value of the integral is well-approximated by the value of the Airy function $(n/\bt^2)^{-1/3}{\rm Ai}((n/\bt^2)^{2/3}\eps)$.   Hence we expect that the integral in \eqref{1} is bounded from below  by $c\cdot (n/\bt^2)^{-1/3}$ for some explicit constant $c>0$, as long as the following Airy condition holds:
\begin{equation}\label{cond}
0\leq \eps \leq (n/\bt^2)^{-2/3}\leq 1/2\ . 
\end{equation} Under this condition, we get the lower bound for the matrix coefficient of the form
\begin{equation}\label{mat-c-bound-sph}
|\langle\pi_{\lambda}(a)e_m,
e_n\rangle|^2 =\pi^{-2}|a|^{-2}|I_{m,n}(a)|^2\geq C\cdot n^{-2/3}\beta^{-2/3}\ ,
\end{equation} for some explicit constant $C>0$ depending on $\lm$, but not on $a$, $m$ and $n$.


To justify the above heuristics,  we use the standard asymptotic for the Airy integral in the formula (36.12.11) of \cite{ol}. 
The partial derivative of $S$ with respect to $x$ is $\frac{\partial S}{\partial x}=\frac{1}{1+x^2}-\frac{(1+\eps)^{-2}}{1+x^2(1+\eps)^{-4}\beta^{-4}}$.
There are two critical points for $S$ in the interval $[-1,1]$:
$x_{\pm}=\pm\sqrt{\frac{2\eps+\eps^2}{1-(1+\eps)^{-2}\beta^{-4}}}$.  Assume that $\eps$
tends to 0 and $\beta$ tends to infinity, as $n$ goes to infinity. Then $x_{\pm}\sim\pm \sqrt{2\eps}$, when $n$ is large.
Substituting values of $x_\pm$ into the formula (36.12.11) of \cite{ol}, we get
\begin{itemize}
\item[(1)] $\tilde{S}=\frac12\big(S(\beta,\eps,x_+)+S(\beta,\eps,x_-)\big)$,
\item[(2)]
$f_{\pm}=f(\beta,\eps,x_{\pm})$,
\item[(3)]
$S^{\prime\prime}_{\pm}=
\frac{\partial^2\!S}{\partial x^2}(\beta,\eps,x_{\pm})$,
\item[(4)]
$\Delta=\left(\frac34[S(\beta,\eps,x_-)-S(\beta,\eps,x_+)]\right)^{2/3}$.
\end{itemize} 
For $n$ large, this gives
\begin{itemize}
\item[(1')]
$\tilde{S}\equiv0$ ($S$ is odd in $x$), 
\item[(2')]
$f_+=f_-\sim 1$ ($f$ is even in $x$), 
\item[(3')]
$S_+^{\prime\prime}=
-S_-^{\prime\prime}\sim
-2\sqrt{2\eps}$ ($S_x^{\prime\prime}$ is odd in $x$),
\item[(4')]
$\Delta\sim
\left(\frac32\right)^{2/3}\cdot2\eps$.
\end{itemize} 
Hence we obtain the following asymptotic (see 36.12.11 of \cite{ol}):
\begin{eqnarray*}&&
I_{m,n}(a)=\frac{\Delta^{1/4}\pi\sqrt{2}}{m^{1/3}}\exp\left(im\tilde{S}\right)
\Big\{\left(\frac{f_+}{\sqrt{S_+^{\prime\prime}}}
+\frac{f_-}{\sqrt{S_-^{\prime\prime}}}\right){\rm Ai}\left(-m^{2/3}\Delta\right)
(1\\ &&+\mathcal{O}(1/m))\Big\} -
i
\left(\frac{f_+}{\sqrt{S_+^{\prime\prime}}}-
\frac{f_-}{\sqrt{S_-^{\prime\prime}}}\right)
{\rm Ai}^{\prime}\left(-m^{2/3}\Delta\right)
\left(1+\mathcal{O}(1/m)\right) \ ,
\end{eqnarray*}
which yields
\begin{eqnarray*}
I_{m,n}(a)&=&\frac{2\pi\sqrt{2}}{m^{1/3}}\cdot\frac{\Delta^{1/4}}{\sqrt{S^{\prime\prime}_+}}\,{\rm Ai}\left(-m^{2/3}\Delta\right)
\big(1+O(1/m)\big)\\[0.2cm]
&\sim&2\pi\left(\tfrac32\right)^{1/6}m^{-1/3}
{\rm Ai}\left(-\left(\tfrac32\right)^{2/3}m^{2/3}\cdot2\eps\right)\ .
\end{eqnarray*} 
Note that $\Delta$, $\Delta^{1/4}$ and $\sqrt{S^{\prime\prime}_+}$
have branched values. When  applying the above asymptotic  formula for $I_{m,n}(a)$, we should select the branches such that both $\Delta$ and $\frac{\Delta^{1/4}}{\sqrt{S^{\prime\prime}_+}}$ are real and positive (see \cite{ol}). To guarantee that the argument in
${\rm Ai}\left(-\left(\tfrac32\right)^{2/3}m^{2/3}\cdot2\eps\right)$
is close to $0$ (and hence the value of ${\rm Ai}$ is bounded away from zero), we have to set the condition \eqref{cond}.


We have shown that the following lower bound holds for the matrix coefficient
\begin{proposition}\label{m-c-bound-prop} Let $a\geq 1$, $m, n>0$ be such that the parameters $\beta=(n/m)^\haf$ and $\eps=a/\beta -1$ satisfy the Airy condition \eqref{cond}.  We have then the following lower bound
\begin{equation}\label{mat-c-bound}
|\langle\pi_{\lambda}(a)e_m,
e_n\rangle|^2 \geq C\cdot n^{-2/3}\beta^{-2/3}\ ,
\end{equation} for some explicit constant $C>0$ depending on $\lm$, but not on $a$, $m$ and $n$.
\end{proposition}

\subsubsection{Test function $\al$}\label{test-al} We now want to integrate the matrix coefficient over a ``spherical shell'' inside $G$. Here we construct our test function $\al\in C^\8_0(G//K)$. Let $\al_1$ be a smooth real valued non-negative function such that $supp(\al_1)\subset [0,1]$, $\al_1(t)\leq 1$ and $\al_1(t)\equiv 1$ for $t\in [1/4,3/4]$.  For a given $\eps>0$, define $\al_\eps$ by $\al_\eps(t)=\al_1(\eps\inv (t-1))$. We have $supp(\al_\eps)\subset [1,1+\eps]$. For a given $T>1$, we define $\al(a)=\al_{\eps,T}(a)=\al_\eps(T\inv a)$. The resulting function is a ``smooth characteristic function'' of the interval $[T,T(1+\varepsilon)]$. We have $\al(t)\equiv 1$ for $t\in [T(1+\eps/4),T(1+3\eps/4)]$.  Note that the value of the integral $\int_{A^-}
\al(a)\rho(a)da$ is of order of  $T^2\varepsilon$
when $T$ is large ($\rho(a)\sim a^{-3}/2$ for $a\inv$ large), and this is of the same order as the volume of the ``spherical shell'' $S_{T,\eps}=KA_{[T,T(1+\varepsilon)]}K\subset G$ (for $\eps\leq 1$).

Since the matrix coefficient in Proposition \ref{m-c-bound-prop} has the same order for $a$ in the support of $\al=\al_{\eps,T}$, we have
for some explicit constants $c_1, c>0$, 
\begin{equation*} \int_{A^-}
\left|\langle\pi_{\lambda}(a)e_m,
e_n\rangle\right|^2\al(a)\rho(a)da\geq c_1 \int_{A^-}
\al(a)\rho(a)da\cdot n^{-2/3}T^{-2/3}\geq c \cdot \eps n^{-2/3}T^{4/3}\ ,
\end{equation*} for $n$, $m$, $T$ and $\eps$ satisfying the Airy condition \eqref{cond}.

We arrive at
\begin{corollary}\label{cor-sph} For $a$,  $m$ and $n$ satisfying conditions of Proposition \ref{m-c-bound-prop} and for the function $\al=\al_{\eps,T}$ constructed above, the following bound holds:
\begin{equation}\label{for-em2} \int_{A^-}
\left|\langle\pi_{\lambda}(a)e_m,
e_n\rangle\right|^2\al(a)\rho(a)da\geq  c \cdot \eps n^{-2/3}T^{4/3},
\end{equation} for some explicit constant $c>0$ depending only on $\lm$ and on the choice of function $\al_1$ as above.  

\end{corollary}

\subsection{Geometric side}\label{sel-sect}
In this section, we obtain an upper bound on the left hand side of (\ref{eq}) (we call it a ``geometric bound" although the main ingredient is  Selberg's {\it spectral} method for lattice points counting).

\subsubsection{Lattice points counting}\label{lat-pt-sel} The function  $H_{\al,x_0}(x)=\sum_{\gamma\in\Gamma}
\al(g_0\inv\gamma g)$ from \eqref{fold-sph} naturally appears in the classical problem of lattice points counting in the hyperbolic plane (initiated by J. Delsarte \cite{De}, H. Huber \cite{hu} and A. Selberg \cite{se}). In particular, Selberg developed a spectral approach which gives an asymptotic formula with a remainder.

The lattice points counting problem we are dealing with has roughly the form
$$\#\{\g\in\G\ |\ g_0\inv\g g\in supp(\al)\}\ .$$ Here $g_0$ is fixed,  $g\in G$, and both could be taken in  $\G\sm G$. We will be interested in a pointwise as well as  the averaged over $g\in \G\sm G$ bound on this function. We deduce the necessary bound from the counting of lattice points in big balls (although one can apply directly the method to the function $\al$). The main point is that one has strong bounds for the remainder in the asymptotic formula. Let $\chi_T={\bf 1}_{[1,T]}$ be the characteristic function of the interval $[1,T]$, $T>1$. We view $\chi_T$ as a bi-$K$-invariant function on $G$ via the Cartan $KA^-K$-decomposition. Selberg (see \cite{se}) obtained the following pointwise asymptotic for the value of $H_{\chi_T,g_0}(g)$
\begin{eqnarray}\label{sel} && H_{\chi_T,g_0}(g)=\\ && vol(X)\inv\int_G\chi_T(g')dg'\cdot {\bf 1}+\sum\limits_{\lm_j\in(0,1)}c_j\phi_{\lm_j}(g_0)\bar\phi_{\lm_j}(g)T^{1+\lm_j}+O(T^{4/3})\ .\nonumber\end{eqnarray}
Here the summation is over a (finite) orthonormal set of  exceptional eigenfunctions on $X/K$ (i.e., eigenfunctions with an eigenvalue $0<\mu_j=(1-\lm_j^2)/4< 1/4$), coefficients $c_j$ are explicit and depend only on $\lm_j$. We note that the main term satisfies $\int_G\chi_T(g')dg'\sim T^2$ and gives the contribution from the $L^2$-normalized constant eigenfunction $\phi_1\equiv vol(X)^{-\haf}$ on $X$ (so the sum above is in fact over $\lm_j\in(1/2,1]$). Note also that functions $\phi_{\lm_j}$ are   bounded for $\G$ co-compact and have some explicit growth in the cusp for non-uniform $\G$. It is expected that the remainder above is of the order of $T$ (and then it would be optimal), but no progress has been made on this conjecture for any $\G$ (see a discussion in \cite{pr}). R. Hill and L. Parnowski in \cite{hp} showed that the average form of this conjecture  holds. Namely, we have, for some explicit constants $c_j'$ depending only on $s_j$,  
$$||H_{\chi_T,g_0}(g)-vol(X)\inv\int_G\chi_T(g')dg'\cdot{\bf 1}-\sum\limits_{\lm_j\in(0,1)}c'_j\phi_{\lm_j}(g_0)\bar\phi_{\lm_j}(g)T^{1+\lm_j}||_{L^2(X)}=O(T)\ .$$ Hence we can represent the function $H_{\chi_T,g_0}(g)\in L^2(X)$ as the sum
\begin{eqnarray}\label{H-rem}
&&H_{\chi_T,g_0}(g)=\\ && vol(X)\inv\int_G\chi_T(g')dg'\cdot {\bf 1}+\sum\limits_{\lm_j\in(0,1)}c'_j\phi_{\lm_j}(g_0)\bar\phi_{\lm_j}(g)T^{1+\lm_j} +R_{g_0,T}(g)\nonumber
\end{eqnarray} with the function $R_{g_0,T}\in L^2(X)$ satisfying $||R_{g_0,T}||_{L^2(X)}\leq C T$ for some constant $C>0$. We now apply this decomposition to the function $H_{\al,g_0}(g)$. Let $\al$ be a smooth characteristic function of the interval $[T,T(1+\varepsilon)]$ as in Corollary \ref{cor-sph}. We have $0\leq \al(g)\leq \chi_{T(1+\varepsilon)}(g)-\chi_T(g)$ and  hence  $H_{\al,g_0}(g)\leq H_{\chi_{T(1+\varepsilon)},g_0}(g)-H_{\chi_T,g_0}(g)$ for all $g\in G$. This together with \eqref{H-rem} implies that
\begin{eqnarray}\label{H-al-rem}
&&0\leq H_{\al,g_0}(g)\leq vol(X)\inv\int_G[\chi_{T(1+\varepsilon)}(g')-\chi_T(g')] dg'\cdot {\bf 1}+\\ &&\sum\limits_{\lm_j\in(0,1)}c'_j\phi_{\lm_j}(g_0)\bar\phi_{\lm_j}(g)T^{1+\lm_j}[(1+\varepsilon)^{1+\lm_j}-1] +R'_{g_0,T}(g)\ , \nonumber
\end{eqnarray}
with $||R'_{g_0,T}||_{L^2(X)}\leq C' T$. Note that for a co-compact $\G$, functions $\phi_{\lm_j}$ are bounded and hence the first term in the decomposition majorates all terms in the sum over the exceptional spectrum for all $T$ and $\varepsilon>0$. (This is not true, of course, for the remainder $R'_{g_0,T}$ which is bigger than the main term if $\varepsilon$ is too small.) Also note that the leading term is of the order of $vol(X)\inv\int_G\al(g')dg'\sim \eps T^2$, i.e., we have 
\begin{equation}\label{lead-bd}
\int_G[\chi_{T(1+\varepsilon)}(g')-\chi_T(g')] dg'\leq A\cdot\int_G\al(g')dg' \sim \eps T^2
\end{equation} for some $A>0$ independent of $T$ and $\varepsilon>0$. 

\subsubsection{$L^4$-norm of  $K$-types} The following theorem gives a uniform upper bound on the $L^2$-normalized  $K$-types in the principal series automorphic representation (see Theorem 2.6 of \cite{br}):
\begin{theorem}\label{L4-norm} Let $(\pi,\nu)$ be an automorphic representation of the principal series and $e_m\in V_\pi$ a norm one vector of  $K$-type $m$. 
 There exists a constant $C=C(\pi,\nu)\geq 1$ such that $$\|\phi_{e_m}\|_{L^4(X)}\leqslant C$$
holds for all $m\in\mathbb{Z}$.
\end{theorem}
We note that the proof of this theorem is also based on the Airy type phenomenon for triple product matrix coefficients and on certain (highly non-trivial) period identities arising form Gelfand pairs. 
 
\subsubsection{Proof of Theorem \ref{sph-thm}}\label{proof-thm-sph} We assume that $\G$ is co-compact, and hence functions $\phi_{\lm_j}$  in \eqref{H-al-rem} are bounded. From \eqref{H-al-rem} we deduce that (taking into account that both functions $|\phi_{e_m}(x)|^2$ and $ H_{\al,x_0}$ are non-negative)
\begin{eqnarray}\label{geom}\nonumber
&&\int_{\Gamma\backslash G}
|\phi_{e_m}(x)|^2 H_{\al,x_0}(x)dx \\ \nonumber&&\leq\Big|vol(X)\inv\int_G[\chi_{T(1+\varepsilon)}(g')-\chi_T(g')]dg'\cdot {\langle |\phi_{e_m}|^2,\bf 1\rangle}+ \\ && \sum\limits_{\lm_j\in(0,1)}c'_j, \phi_{\lm_j}(g_0)\langle |\phi_{e_m}|^2,\bar\phi_{\lm_j}\rangle T^{1+\lm_j}[(1+\varepsilon)^{1+\lm_j}-1]+  \langle |\phi_{e_m}|^2, R'_{g_0,T}\rangle\Big|\ .\nonumber
\end{eqnarray} The first term on the right could be bounded by $vol(X)\inv A\cdot\int_G\al(g')dg'$ for some $A>0$ by \eqref{lead-bd}. All functions $\phi_{\lm_j}$ are bounded (since $\G$ is co-compact) and $s_j< 1$. Hence we have that the sum over the exceptional spectrum is also bounded by $A'\cdot\int_G\al(g')dg'$ for some $A'>0$. Hence the sum over the spectrum below $1/4$ is of order of $\eps T^2$. We bound the last term  by the Cauchy-Schwartz inequality as  $|\langle |\phi_{e_m}|^2, R'_{g_0,T}\rangle|\leq ||\phi_{e_m}||_{L^4(X)}|| R'_{g_0,T}||_{L^2(X)}\leq D T$ for some $D>0$. Hence we obtain the geometric bound 
\begin{equation}\label{geom-fin}
\int_{\Gamma\backslash G}
|\phi_{e_m}(x)|^2 H_{\al,x_0}(x)dx\leq B\eps T^2+D T\ ,
\end{equation}
with some explicit constants $B, D>0$ depending on $\G$, $x_0$, $(\pi,\nu)$ and the function $\al'$ we used to construct the function $\al$ in Corollary \ref{cor-sph}, but not on $m$ and $T$. This bound together with the lower bound \eqref{for-em2} on the integral of the matrix coefficient and the identity \eqref{eq} implies the bound \eqref{sph-thm-bound} $|c_n|^2\leq C |n|^{6/7}$. Namely, dropping all but one term in the sum \eqref{eq}
\begin{equation*}\sum\limits_n|c_n|^2\cdot\int_{A^-}
\left|\langle\pi_{\lambda}(a)e_m,
e_n\rangle\right|^2\al(a)\rho(a)da\ 
\end{equation*}
and using lower bound \eqref{for-em2} and \eqref{geom-fin},  we have for $T=\bt=(n/m)^\haf>1$ and  constants $B, D>0$ from  above,
\begin{equation}\label{inq1}|c_n|^2\cdot \eps n^{-2/3}T^{4/3}\leq B\eps T^2+ D T\ ,
\end{equation}
or $|c_n|^2\leq  n^{2/3}(BT^{2/3}+ D \eps\inv T^{-1/3})$. 
We also have to satisfy the Airy condition $0\leq\eps\leq n^{-2/3}T^{4/3}\leq 1$ from Corollary \ref{cor-sph}. Balancing two terms (with $T=n^{2/7}$ and hence $m=n^{3/7}$ or $m=T^{3/2}$), we obtain the bound desired. 

To obtain a bound for the sum of coefficients $|c_n|^2$, we note that for a given $v=e_m$, integrated matrix coefficients $|\langle \pi(a)e_m,e_{n'}\rangle|^2$ in the identity \eqref{eq} satisfy the bound \eqref{for-em2} (maybe with a smaller constant $c>0$) as long as the Airy condition in the integral $I_{m,n}(a)$ in \eqref{I-mna}  is satisfied. For a given $m$, the condition reads $|1-n'\bt^{-2}m\inv(1+2\eps)|\leq\eps$. For a given $\bt=T>1$ as above, we have $m= T^{3/2}$, $n=T^{7/2}$ and hence  $\eps\leq n^{-2/3}\bt^{4/3}=T^{-1}$. This implies that the range of $n'$ in \eqref{eq} amplified by a given $v=e_m$ satisfy $|1-n'T^{-7/2}|\leq \eps\leq T^{-1}$ or $|T^{7/2}-n'|\leq T^{5/2}$. Setting $S=T^{7/2}$ we arrive at $|S-n'|\leq S^{5/7}$ as in  \eqref{sph-thm-bound-aver}.

\subsubsection{$L^2$-theory.}\label{L2-th-sect} The above argument is quite flexible and in fact does not require the highly non-trivial result on $L^4$ norms of $K$-types in order to obtain some polynomial saving over the ``trivial'' bound $|c_n|\ll n^\haf$. 

The basic identity \eqref{eq} holds for any vector $v\in V_\pi$. In order to use it to bound   coefficients $c_n$ (as $|n|\to\8$), we need an amplification  
of the matrix coefficient $\langle\pi_\lm(g)v,e_n\rangle$ for $g$ in the support of some function $\al$ and the geometric bound on the counting function $H_{\al,x_0}$. We achieved the amplification for $v=e_m$ using the Airy phenomenon, but there might be other ways for  different vectors $v$ (see Section \ref{non-stand-sph} below). Assuming that the bound in Corollary~\ref{cor-sph} holds, we want a {\it pointwise} bound on $H_{\al,x_0}(g)$  for all $g\in \G\sm G$. Let us assume that we have a very coarse Selberg type bound for the lattice points counting problem in big balls similar to the bound \eqref{sel}: 
\begin{eqnarray}\label{H-rem-gen}
H_{\chi_T,g_0}(g)= vol(X)\inv\int_G\chi_T(g')dg'\cdot {\bf 1} +R_{g_0,T}(g)\nonumber
\end{eqnarray} with the remainder satisfying a {\it pointwise} bound $|R_{g_0,T}(g)|\ll T^{2-\s}$ for some $\s>0$ (note that $\int_G\chi_T(g)dg\sim T^2$). Note that again the contribution over the exceptional spectrum is bounded by the volume $\eps T^2$ of the shell and hence we can assume that the origin of the exponent in the  remainder is not connected with the existence of the exceptional spectrum. Bounding the term $\langle |\phi_{e_m}|^2, R_{g_0,T}\rangle$ in \eqref{geom}\ by the supremum  norm of $R_{g_0,T}$,  our argument above  implies that 
\begin{eqnarray*}\label{geom-L2}&&
\int_{\Gamma\backslash G}
|\phi_{e_m}(x)|^2 H_{\al,x_0}(x)dx \\ && \leq vol(X)\inv\int_G[\chi_{T(1+\varepsilon)}(g')-\chi_T(g')]dg'\cdot {\langle |\phi_{e_m}|^2,\bf 1\rangle}+  \langle |\phi_{e_m}|^2, R'_{g_0,T}\rangle\\ && \leq B_1T^2\eps+D_1T^{2-\s} .\nonumber
\end{eqnarray*}

Hence in this case, the analog of \eqref{inq1} reads 
\begin{equation}\label{inq1-L2} |c_n|^2\cdot \eps n^{-2/3}T^{4/3}\leq B_1\eps T^2+ D_1 T^{2-\s}\ .
\end{equation} It is easy to see that for any $\s>0$, this implies a non-trivial bound $|c_n|^2\leq C' n^{1-\dl(\s)}$ with $\dl(\s)>0$ (to be precise: $\dl(\s)=\s/(4+3\s)$; all this under the Airy condition on $\eps$, $n$ and $T$ in Corollary \ref{cor-sph}). In particular, Selberg's bound \eqref{sel} ($\s=2/3$) gives $|c_n|^2\leq C' n^{8/9}$. 


As in Section \ref{proof-thm-sph}, we  get a bound for the sum of coefficients $|c_n|^2$ over a short interval. For this purpose, we  consider those $n^{\prime}$ for which the matrix coefficient $|\langle\pi(a)e_m,e_{n^{\prime}}\rangle|^2$ satisfy the bound \eqref{mat-c-bound}. For a given $\beta=T\gg 1$, we have from the above argument that $m=aT$, $n=b T^3$ and $\eps=m^{-2/3}=c T^{-2/3}$ for some $a,b,c>0$. With $n$ replaced by $n^{\prime}$ such that $|1-n^{\prime}/n|\leq \eps$, the Airy condition is still satisfied. Hence, for such $n^{\prime}$ we have 
$$\eps n^{-2/3}T^{4/3} \sum_{|1-n^{\prime}/n|\leq \eps}|c_{n^{\prime}}|^2 \leq B_1\eps T^2+ D_1T^{2-\sigma}.$$ 
With $n$, $\eps$ as above, the condition $|1-n^{\prime}/n|\leq \eps$ reads $|1-n^{\prime}c^{-1}T^{-3}|\leq e T^{-2/3}$ which is equivalent to $|n^{\prime}-cT^3|\leq e'T^{7/3}$. Let $S=cT^3$, then $|n^{\prime}-S|\leq e'' S^{7/9}$ for some $e''>0$. As in Section \ref{proof-thm-sph} we arrive at $\sum_{|n-S|\leq  S^{7/9}}|c_n|^2\leq C^{\prime}S^{8/9}$.

\subsubsection{Non-spherical counting}\label{non-stand-sph} Here we briefly discuss another choice of a test vector in \eqref{eq} which leads to a ``non-classical'' counting problem. While this would not give a better exponent in the bound \eqref{sph-thm-bound}, we feel that it shows certain flexibility of the method. In particular, we reduce the problem to counting in domains defined with respect to Iwasawa decomposition as opposed to Cartan decomposition we used before.
This time the amplification is achieved by the use of a test vector with a monomial {\it real} singularity (i.e., closeness of a monomial singularity of the amplitude to a non-degenerate critical point of the phase functions). This is    another classical phenomenon in the theory of oscillating integrals which is different from the Airy type singularity (i.e., closeness of two non-degenerate critical points of the phase function).  

Consider a vector $v_{\xi,\s}\in V_\pi$ given in the line model by  $v_{\xi,\s}(x)=|x|^{-\s}e^{i\xi x}\chi(x)$, where $\haf>\s>0$, $
\xi\in\RR$,  and $\chi\in C^\8_0(\RR)$ is a smooth compactly supported function which is identically one in a neighborhood of zero (e.g., $supp(\chi)\subset [-1,1]$ and $\chi|_{[-\haf,\haf]}\equiv 1$). Note that the vector $v_{\xi,\s}$ is in $L^2$, but it is not a smooth vector.

We have  the  following well-known asymptotic for the oscillatory integral  defined by the matrix coefficient $ \langle\pi([a][n_t]_-)v_{\xi,\s},e_n\rangle=\langle\pi([n_t]_-)v_\xi,\pi([a]\inv)e_n\rangle$. Here $|n|\to\8$ is a parameter, $\xi\in\RR$ is an auxiliary parameter, $[a]= {\rm diag}(a,a\inv)$ and $[n_t]_-=\left(\begin{smallmatrix}1&0\\t&1
\end{smallmatrix}\right)$ where $t$ is in a {\it fixed} small neighborhood of zero (e.g., $t\in[-1/4,1/4]$).

We have $\pi_\lm([n_t]_-)v(x)=|1-tx|^{\lm-1}v(\frac{x}{1-tx})$ and hence $\langle\pi([n_t]_-)v_\xi,\pi([a]\inv)e_n\rangle=|a|^{1-\lm}m_{\xi,\s,n}(a,t)$ with
\begin{eqnarray}\label{matrix-c-vxi} m_{\xi,\s,n}(a,t)=
\int\limits_\RR|1-tx|^{\lm-1}\left|\frac{x}{1-tx}\right|^{-\s}\chi\left(\frac{x}{1-tx}\right)e^{i\xi(\frac{x}{1-tx})}e^{-in\arctan(a^{2}x)}dx\ .\nonumber\\
\end{eqnarray}
Choosing the support of $\chi$ and the size of $t$ small enough, we can assume that $1-tx\not=0$. We also assume that  $\lm=\lm_0$ is fixed. Hence the above integral is an oscillating integral with the {\it singular} amplitude $|x|^{-\s}$. We write the integral \eqref{matrix-c-vxi} in the standard form  $m_{\xi,n}(a,t)=\int |x|^{-\s}e^{inS(\xi,n,a,t)}f_{t}(x)dx$. We will set  parameters $a$ and $\xi$ so that the phase function has the only critical point near $x=0$. We have $S(\xi,n,a,t)=    \frac{\xi}{n}(x+tx^2)-a^2(x+\frac{a^4}{6}x^3)+{\rm smaller\ order\ terms\ in\ } x\ {\rm and}\ a$. By setting $\xi=a^2n$, we see that the integral is well approximated by the standard integral $\int |x|^{-\s}e^{ia^2n tx^2}dx$ (see \cite{mo}) and is of the order of $|a^2n|^{\frac{\s}{2}-\haf}$ (for $t$ bounded away from zero). 
It is very easy to justify the above heuristic and to show formally that for $ 1>a=\bt(1+\eps)\gg n^{-\haf} $, $\xi=a^2n>1$ and $\haf-\s=\dl>0$ arbitrary small, we have
\begin{equation}\label{mat-c-bound-xi}
|a\cdot m_{\xi,\s,n}(a,t)|^2=|\langle\pi_{\lambda}([a][n_t]_-)v_{\xi,\s},
e_n\rangle|^2 \geq C_\dl\cdot n^{-1/2-\dl}\beta^{1-2\dl}\ ,
\end{equation} for some  constant $C_\dl>0$ depending on $\dl>0$ and $\lm$, but not on $a$ and $n$. This bound is valid under the condition $0<\eps<(\bt^2n)^{-\haf}$ (this ensures that the critical point of the phase is at most at a distance $(\bt^2n)^{-\haf}$ from the zero, which is the {\it real} singularity of the amplitude $|x|^{-\s}$). 

We now set up the counting problem arguing as in the spherical case, but this time choosing the left $K$-invariant function $\al_{T,\eps}\in C^\8_0(G)$ supported in a ``shell'' $B_{T,\eps}$ in $KAN_-$-coordinates with the $N_-$ part in a fixed  bounded set in $N^-$ and $[a]\in A$ such that $a\in [T\inv, T\inv(1+\eps)]$.   Note that the Haar measure in $k[a][n_t]_-$-coordinates is given by $a^{-2}dk\frac{da}{a}dt$ and hence the volume of this shell is of order  of $T^2\eps$. We now assume that we have an effective (that is, with the main term given by the volume and a remainder with a power saving) asymptotic for the lattice point counting in balls $B_T=KA^-_TN^0_-\subset G$ given in $KAN_-$ decomposition by taking a compact set $N^0_-\subset N_-$ and $[a]\in A$ with $a\in [T\inv, 1]$. Such a counting problem appeared under the name ``counting in sectors'' (see \cite{go}, \cite{go}, \cite{gos}). Effective lattice point counting  would imply that we have a pointwise bound for the counting function $H_{\al_{T,\eps}}$ of the type
\begin{eqnarray*}\label{non-sph-H-al}
H_{\al_{T,\eps}}(x)dx\leq  BT^2\eps+DT^{2-\vartheta} ,\nonumber
\end{eqnarray*} for some  constants $B, D>0$ and some $\vartheta>0$. Arguing as in the spherical case, we would arrive at the bound 
\begin{equation}\label{inq1-L2-non-stand} |c_n|^2\cdot \eps n^{-\haf-\dl}T^{1-2\dl}\leq B\eps T^2+ DT^{2-\vartheta}\ ,
\end{equation} for {\it any} $\dl>0$, some {\it fixed} $\vartheta>0$ and under the condition $0<\eps<(T^{-2}n)^{-\haf}$. It is easy to see that such a bound would imply a non-trivial bound $|c_n|^2\leq C|n|^{1-\vartheta'}$ for some explicit 
$\vartheta'>0$.

\subsubsection{Non-uniform lattices}\label{non-comp} We now consider non-uniform lattices $\G\subset G$.  In this case $X$  has finitely many cusps. In this section, we slightly alter our choice of the test function $\al$ from  Section \ref{test-al} and apply Selberg's pointwise lattice counting asymptotic to this function to show how to derive a non-trivial bound on spherical periods.  The standard observation is that the contribution from cusps is negligible and hence the same argument works for non-uniform lattices as well.

For simplicity, we assume that $\infty$ is the unique  cusp of $\G$which is reduced at infinity, i.e., $\Gamma_{\infty}:=\Gamma\cap N=\left\{\begin{pmatrix}1&n\\0&1\end{pmatrix}\middle|n\in\mathbb{Z}\right\}$. The main difference with the co-compact case we treated before is that Selberg type lattice points counting  asymptotic \eqref{sel} have terms growing in the cusp. These are controlled by Eisenstein series and residual spectrum. Exceptional eigenfunctions (i.e., those $\phi_{\lm_j}$ with $\lm_j\in(0,1)$ in \eqref{sel}) satisfy the following growth condition: $|\phi_{\lm_j}(z)|\leq a_j\cdot y^{\haf-\haf\lm_j}$, $z=x+iy$, for some $a_j>0$. Let $\s=\s_\G=\max\{\lm_1, 4/3\}$ be the ``counting'' spectral gap of $\G$. The corresponding spherical pointwise counting asymptotic \eqref{sel}  is (in a coarse form) the following statement
\begin{eqnarray}\label{sel-non-uni}  H_{\chi_T,g_0}(g)= vol(X)\inv\int_G\chi_T(g')dg'\cdot {\bf 1}+O_{g_0}(y^{\haf}T^{1+\s})\ .\end{eqnarray}

For $A\geq 1$ to be chosen later, denote by $\Omega_A$ the neighborhood of the cusp $\Omega_A=\{x+iy\,:\,|x|\leq 1/2,~y\geq A\}\subseteq\mathbb{H}\simeq G/K$ in the fundamental domain $\Om$ of $\Gamma$. We have ${\rm vol}(\Omega_A)=A^{-1}$. Let $\tau_{A}$ be a smooth characteristic function of $\Omega_A$. Define $\beta_{T,A}\in C^{\infty}_c(G)$ to be a non-negative right $K$-invariant smooth function given for $g\in G$, by
$$\beta_{T,A}(g)=\chi_{T}(g)-\sum_{\gamma\in\Gamma}\tau_{_A}(\gamma^{-1} g),\quad   $$ if this quantity is non-negative and $\beta_{T,A}(g)=0$ otherwise. 
For any $g\in G$, there exists exactly one $\gamma\in\Gamma$ such that  $\gamma^{-1}g$ lies in $\Omega$. Hence, for any given $g\in G$, the value of $\tau_{_A}(\gamma^{-1}g)$ does not vanish for at most one $\gamma$. We have then $\beta_{T,A}(g)=\chi_{_T}(g)-\tau_{_A}(\gamma_g^{-1}g)$. When $\gamma_g$ does not exist, we simply have  $\beta_{T,A}(g)=\chi_{_T}(g)$. The support of $\beta_{T,A}$ contains the  set (maybe for a slightly smaller value of $A$) $$D_{T,A}=B_T\smallsetminus\left(\left(\bigcup_{\gamma\in\Gamma}\gamma\Omega_A\right)\bigcap B_T\right).$$ The set $D_{T,A}$ could be described as the ball $B_T\subset G$ with all preimages of $\Om_A$ removed.  
Note that the $\gamma\Omega_A$'s are disjoint from each other for different $\gamma$'s. Clearly the function $H_{\bt_{T,A},g_0}(g)=\sum_\g \bt_{T,A}(g_0\inv\g g)$ is supported 
in the compact part of the fundamental domain $\Om\sm\Om_A\subset \G\sm G$ and is bounded. Moreover, the asymptotic \eqref{sel-non-uni} implies that
\begin{equation}\label{cut-sel} H_{\bt_{T,A},g_0}(g)= vol(X)\inv vol(D_{T,A}) {\bf 1}+O_{g_0}(A^{\haf}T^{1+\s})\ , \end{equation} for $g\in \Om\sm\Om_A$ (where the averaging over $\G$ is not affected by our change of the test function) and zero otherwise. We claim that for a large enough fixed $A$, the set $D_{T,A}$ is still a large part of the ball $B_T$ and hence the first term in the asymptotic 
\eqref{cut-sel} dominates the second term. Namely we claim that for $A$ which is a large enough fixed number (depending on the geometry of $\Om$), we have $vol(D_{T,A})\sim C_A T^2$ for some $C_A>0$. In fact, this again follows from the asymptotic \eqref{sel-non-uni} which implies that the total volume of the set $\left(\bigcup_{\gamma\in\Gamma}\gamma\Omega_A\right)\bigcap B_T$ is bounded by $A^{-\haf} T^{1+\s}$. This gives a pointwise bound for the averaging of the characteristic function for the shell $B_{T(1+\eps)}\sm B_T$ with preimages of $\Om_A$ removed, and hence we can apply the argument in Section~\ref{L2-th-sect} to the function $\bt_{T,A}$ instead of the function  $\al_T$ we used in Section \ref{test-al}.


\section{Hyperbolic periods}\label{hp}
In this section we deal with hyperbolic periods, namely, those periods arising from the diagonal subgroup $A$.  

\subsection{Closed geodesics}  We denote by $G^+$ the neutral  connected component of $G=\PGLR$, and by $G^-$ the other connected component. In what follows, it is crucial that we work with both components.  Let $Y$ be the Riemann surface as before.  The uniformization  theorem implies that there exists a lattice $\G\subset G^+$  such that $Y\simeq\G\setminus G^+/K$.  Let $X=\G\setminus G$ be the automorphic space and  $X_\pm=\G\setminus G^\pm$ are two of its connected components. The Riemann surface $Y_-=X_-/K$ is naturally identified with the orientation reversed copy of $Y$, and $X_\pm$ could be viewed as the spherical (co-)tangent bundle of $Y_\pm$. 
Let $A=\{{\rm diag}(a,b)\}/\{{\rm diag}(a,a)\}\subset G$ be the full (disconnected) diagonal subgroup, and denote by $A_\pm=A\cap G^\pm$ its connected components. A closed geodesic $l\subset Y$ corresponds to an $A_+$-orbit  in $X_+$, which we denote by the same letter $l$. We fix an $A_+$-invariant measure $dl$ on $l$ (e.g., consistent with the Riemannian length of $l$). Let $\dl=\left(
                                          \begin{smallmatrix}
-1  & \\
                                   &\ \ 1 \\
\end{smallmatrix}
                                        \right)\in \GLR$,
and we denote by the same letter the corresponding element in $G$. We have $A=A_+\cdot\langle\dl\rangle$.

We now consider closed $A$-orbits. If $l\subset X_+$ is a closed $A_+$-orbit, then $\OO=l\cdot A=l\cup l\cdot\dl\subset X$ is a closed $A$-orbit consisting of (possibly) two connected components $l_\pm\subset X_\pm$ interchanged by $\dl$.
We have then the {\it pointwise} stabilizer  subgroup $A_\OO={\rm Stab}_A  \OO={\rm Stab}_{A_+}  l=\langle [a]_q\rangle \simeq \ZZ$,
 where $ [a]_q={\rm diag}(e^q,e^{-q})\in A_+$, $q>0$, is an element which is conjugated to a primitive hyperbolic element $\g_l\in \G$ corresponding to the closed geodesic $l$ on $Y$ (i.e., $g_l\inv \g_l g_l= [a]_q$ for an appropriate $g_l\in G$).

 \subsubsection{Characters} We consider (unitary) characters $\chi:A\to\CC^\times$ of the disconnected group $A$. Any such unitary character is given by $\chi_{s,\epsilon}\left([a]\right)=|a|^s\cdot {\rm sign}(a)^\epsilon$, where $[ a]={\rm diag}(a,|a|\inv)\in A$, $s\in i\RR$, and $\epsilon=0, \ 1$.  A character  is called  {\it even} if $\epsilon=0$  and {\it odd} if $\epsilon=1$ (note that  $\chi_{s,\epsilon}(\dl)=(-1)^\epsilon$).

 We next consider characters which are trivial on the subgroup $A_\OO$.
These are given by $\chi_{s_n,\epsilon}$ with $s_n=2\pi in/q$, $n\in \ZZ$, and $\epsilon\in\{0,1\}$.

We now construct $A$-equivariant functionals on the space of smooth functions on $X$.  Let $\OO\subset X$ be a closed $A$-orbit as above. We fix a point $x_0=\Gamma\cdot g_0\in \OO\subset X$. We denote by $\G_\OO=g_0\inv \G g_0\cap A$ the corresponding elementary subgroup in $\G$.  We associate to any  character $\chi=\chi_{s_n,\epsilon}$ which is trivial on $A_\OO$ the function on $\OO$ given by $\chi^\cdot(l_0\cdot [a])=\chi([a])$.  This gives rise to an  $A$-equivariant functional $d^{aut}_{\chi}:C^\8(X)\to\CC$  given by $d^{aut}_\chi(f)=\int_{\OO}f(o)\bar\chi^\cdot(o)do$. The functional $d^{aut}_{\chi}$ is $\chi$-equivariant with respect to the right action of $A$: $d_\chi(R(\bar a)f)=\chi(\bar a)f$.

Functions $\{\chi_{s_n,\epsilon}^\cdot\}$ give rise to an orthonormal basis of $L^2(\OO,do)$ where $do$ is the $A$-invariant measure on $\OO$ such that ${\rm vol}(\OO)=1$. For a smooth function $f\in C^\8(X)$, we have the decomposition $\phi|{_{\OO}}=\sum_{n,\epsilon}d^{aut}_{\chi_{n,\epsilon}}(\phi)\chi_{n,\epsilon}^\cdot$, $a_n\in\CC$ which is understood in $L^2$-sense. 

We now consider functionals as above restricted to (irreducible) automorphic representations. The space of such functionals is at most one-dimensional.  Namely,
we have $\dim {\rm Hom}_A(V_{\pi},\chi_{n,\epsilon})\leq 1$ for any irreducible unitary representation (automorphic or not) $(\pi, V_\pi)$ of $\PGLR$ (see \cite{bu}). Note that for the subgroup $A^+$ the space is two-dimensional, and this is what forces us to work with the disconnected group $A$. 

For simplicity, we assume that  the automorphic representation $(\pi, V_\pi,\nu)$ is even representation of the principal series (as an abstract representation of $\PGLR$) and that  the element $\dl\in A$ acts trivially on the orbit $\OO$. This means that there is an orientation reversing automorphism of $Y$ which reverses direction on the geodesic $l$ (i.e., that the corresponding primitive hyperbolic element $\g_l$ is $\G$-conjugate to its inverse) or, equivalently, that $\G_\OO=\langle \g_l, \g_\dl\rangle$ with $\g_\dl$ conjugated to $\dl\in A$.
These two assumptions imply that ``even'' functions $\chi_n^\cdot=\chi_{n,0}^\cdot$ form an orthonormal basis for  $L^2(\OO,do)$. The general case could be treated analogously (see \cite{re1}).  

\subsubsection{Invariant functionals}
For a character $\chi$ as above (i.e., trivial on $A_\OO$), we have the  following functional on an automorphic representation  $(\pi,V_{\pi},\nu)$
\begin{align}\label{aut-hyp-funct}
	d^{\rm aut}_{\chi}:\,V_{\pi}\rightarrow\mathbb{C},\quad v\mapsto\int_{\OO}\phi_v(o)
\overline{\chi}^\cdot(o)do\ .
\end{align}

For a general (even) unitary character $\chi_s=\chi_{is,0}$ of $A$, $s\in \RR$, we  define an equivariant functional by the following integral in the linear model of a representation $\pi\simeq\pi_\lm$   of the principal series (automorphic or not):
\begin{align}\label{mod-hyp-funct}d^{\rm mod}_{s}:\,V_{\lambda}\rightarrow\mathbb{C},\quad v\mapsto\tfrac{1}{\pi}\int_{\mathbb{R}}v(x)|x|^{-\frac{1+is+\lambda}{2}}dx\ .
\end{align}
For a given character $\chi_n$, both functionals   $d^{\rm aut}_n=d^{\rm aut}_{\chi_{s_n}}$ and  $d^{\rm mod}_n=d^{\rm mod}_{\chi_{s_n}}$ belong to the same on-dimensional vector space ${\rm Hom}_A(V_{\pi},\chi_{s_n})$, and hence one-dimensionality implies the existence of a  scalar $c_n=c(s_n,\nu)\in\CC$ such that 
\begin{align}\label{c-hyp}
d^{\rm aut}_{{n}}=c_n\cdot d^{\rm mod}_{{n}} \ .
\end{align}

 The Plancherel formula gives for a vector $v\in V_\pi$,  
\begin{align}\label{Planch-hyp}
||\phi_v|_\OO||^2_{L^2(\OO)}=\sum_n\left|d^{\rm aut}_n(\phi_v)\right|^2=\sum_n|c_n|^2\cdot\left|d^{\rm mod}_n(v)\right|^2\ .
\end{align}
The group $G$ acts on this identity by acting on a vector $v$.  We have 
\begin{equation}\label{eq-4}
||\phi_{\pi(g)v}|_\OO||^2_{L^2(\OO)}=\sum_n|c_n|^2\cdot\left|d^{\rm mod}_n\big(\pi_{\lambda}(g)v\big)\right|^2\ .
\end{equation}

\subsubsection{Airy asymptotic}
Clearly $A$ acts trivially on the  identity \eqref{eq-4}. We will apply this identity to a $K$-type vector $v=e_m$. Hence we can choose $g\in N\simeq A\sm G/K$. We are led to the study of (generalized) matrix coefficients 
\begin{equation}\label{hyp-mat-c}
m_{m,s}(n)=\langle \pi_\lm(n)e_m,d^{\rm mod}_s\rangle=d^{\rm mod}_s(\pi_\lm(n)e_m)\ .
\end{equation} We have $m_{m,s}(n)=\langle \pi_\lm(n)e_m,d^{\rm mod}_s\rangle=\langle e_m,\pi_\lm(n\inv)d^{\rm mod}_s\rangle$.   This gives for $n=n_t=\begin{pmatrix}1&t\\0&1
\end{pmatrix}$, $t\in\RR$, 
\begin{equation}\label{int-hyp-mat-c}
m_{m,s}(t) =\tfrac{1}{\pi}\int_{\mathbb{R}}
e^{im\arctan(x)}\big(1+x^2\big)^{\frac{\lambda-1}{2}}|x+t|^{-\frac{1+is+\lambda}{2}}dx\ . \end{equation}

We write this integral in the standard form $m_{m,s}(t)=\tfrac{1}{\pi}\int_{\mathbb{R}}e^{imS(s,t;x)}f(t,\lm,t; x)dx$, where the phase function is given by 
$S(s,t;x)=\arctan(x)-\frac{s}{2m}\ln|x+t|=\arctan(x)-\frac{s}{2m}(\ln|t|+\ln|1+x/t|)$ and the amplitude function is $f(t,\lm,t;x)=\big(1+x^2\big)^{\frac{\lambda-1}{2}}|x+t|^{-\frac{1+\lambda}{2}}$.  We will be looking for critical points of the phase near $x=0$ for $t$ which is large. We introduce new parameters $T$ and $\eps$ such that $t=T(1+\eps)$ with
$s=2mT$ and $1>\eps>0$. We have then 
\begin{eqnarray}\label{hyp-phase}\nonumber &&
S(s,t;x)=x-\frac{1}{3}x^3-\frac{s}{2mT}(1-\eps)\left(x+\frac{x^2}{2T(1+\eps)}\right)+\\ &&{\it\ smaller\ order\ terms\ in\ } x,\ \eps\ {\rm and}\ T\inv\ .
\end{eqnarray} Hence the phase functions is essentially of the form 
$S(s,t;x)=\eps x+ \frac{x^2}{2T}- \frac{1}{3}x^3$. The value of the amplitude function is of order of $T^{-\haf}$. For $T\gg m^{1/3}$, this is Airy type critical point and hence one expects that 
the matrix coefficient $m_{m,s}(T(1+\eps))$ is asymptotic to  $m^{-1/3}T^{-\haf}{\rm Ai}(m^{2/3}\eps)\approx {s}^{-1/3}T^{-1/6}$ for $0<\eps\leq s^{-2/3}T^{2/3}$. 
Hence we expect that 
\begin{equation}\label{hyp-matr-bd}|m_{m,s}(T(1+\eps))|^2=|\langle \pi_\lm(n_{T(1+\eps)})e_m,d^{\rm mod}_s\rangle|^2\geq C s^{-2/3}T^{-1/3}
\end{equation} for some explicit constant $C>0$, as long as the Airy condition $0<\eps\leq s^{-2/3}T^{2/3}$ is satisfied (here $2mT=s$).
The justification for the bound \eqref{hyp-matr-bd} is easy to obtain from \cite{ol} as in the spherical case.

\subsubsection{Counting}
Let $\al\in C^\8(A\sm G/K)$ be a smooth non-negative  function on $G$ such that it is   right $K$-invariant and left $A$-invariant. We will assume that it is compactly supported on $A\sm G/K$. Note that the left invariant measure on $G$ is given by $d[a]dndk$ in $ANK$ coordinates. 
Integrating \eqref{eq-4} on the left we obtain for a $K$-type vector $v=e_m$, 
\begin{eqnarray*}
&&\int_N\int_\OO\big|\phi_{\pi(n)v}(o)\big)\big|^2\al(n)dodn
=
\int_K\int_N\left(\int_{\OO}\big|
\phi_v(onk)\big|^2dx\right)\al(nk)dodndk =\\[0.2cm]
&&\int_K\int_N\left(\int_{A_\OO\sm A}\big|
\phi_v(g_0ank)\big|^2dx\right)\al(ank)dadndk =\\ &&
\int\limits_{\G_\OO\sm G}\big|
\phi_v(g)\big|^2\al(g_0\inv g)dg
=\int_{X}
\big|
\phi_v(x)\big|^2H_\al(x)dx\ ,
\end{eqnarray*}
where 
$H(x)=\sum_{\gamma\in\G_\OO\sm\Gamma}
\al(g_0\inv\gamma g)$ with $x=\Gamma\cdot g\in X$ and $x_0=\Gamma\cdot g_0\in X$.

Integration on the right side of (\ref{eq-4}) gives 
\begin{equation}\label{hyp-eq}
\sum_n|c_n|^2\int_N\big|\langle\pi_{\lambda}(n)e_m,d^{\rm mod}_{\chi_{s_n}}\rangle\big|^2dn\ .
\end{equation}

To simplify notations, we can assume that $g_0=e$ (i.e., we can conjugate $\G$ inside of $G$ so that $\G_\OO=\G\cap A$).   We now assume that we have an effective counting of lattice points $A\cdot\G\subset A\sm G$  in $K$-invariant balls $B_T\subset A\sm G$ of the form $A\cdot N_TK$ where  $N_T=[0,T]$ under the identification $N_T\subset N\simeq \RR$. Note that the volume of $B_T$ is of the order of $T$. The corresponding counting problem corresponds to counting lattice points on a {\it one-sheeted} hyperboloid and was introduced in \cite{drs} and \cite{em}.  We expect that the number of lattice points in $B_T\subset A\sm G$ satisfies the asymptotic  $\frac{m(l)}{vol(X)}vol(B_T)+ O(T^{1-\vartheta})$ for some $\vartheta>0$. Here $m(l)={\rm length}(l)=vol_{d[a]}(A_l\sm A)$ is the length of the closed geodesic $l$ with respect to the measure $d[a]$ on $A$ such that the Haar measure $dg$ on $G$ is given by $dg=d[a]dndk$ in $ANK$ decomposition. 

Let $\al_{T,\eps}\in C^\8(A\sm G/K)$ be a smooth non-negative characteristic function of the shell $S_{T,\eps}=B_{T(1+\eps)}\sm B_T$. Note that the volume of $S_{T,\eps}$ is of order of $\eps T$. The effective counting on $A\sm G$ would imply that $\sup_{x\in X}|H_{\al_{T,\eps}}(x)|\leq B\cdot \eps T +D T^{1-\vartheta}$ for some $\vartheta>0$ and $B, D>0$. The bound \eqref{hyp-matr-bd} on matrix coefficients then implies that the bound
 \begin{equation}\label{inq1-L2-hyp} |c_n|^2\cdot |s_n|^{-2/3}T^{2/3}\eps\leq B\eps T+ DT^{1-\vartheta}
\end{equation} holds under the condition $0<\eps<|s_n|^{-2/3}T^{2/3}$. This implies a non-trivial bound $|c_n|^2\leq C |s_n|^{1-\vartheta'}$ with $\vartheta'=\vartheta/(2+3\vartheta)>0$ and some $C>0$.

We have proved the following 
\begin{theorem}\label{hyp-thm} There exists a constant $\theta=\theta_{\G,\nu}>0$  depending on $\Gamma$ and $\nu$  such that 
$|c_n|^2\leq C{\cdot}|n|^{1-\theta}$ for some $C>0$. 
\end{theorem}


 \section{Unipotent periods}\label{uni}
In this section we deal with  unipotent periods, namely, those periods arising from closed orbits of the unipotent subgroup $N$. 

\subsection {Integrated Hermitian form}
Let $x_0\in X$ be a point such that the orbit $\OO=x_0N\subset X$ is closed. Namely, $\OO$ is a horocycle in $X$. By conjugating $\G$, we can assume that that $x_0=\G e\in X$ is the image of the identity in $X=\G\sm G$. Moreover, we will assume for simplicity that the lattice $\G$ is reduced at $\8$, that is, it has the standard subgroup fixing the cusp at infinity: 
$\G_\8=\G\cap N=\langle \g_\8\rangle\simeq\ZZ$, where $\g_\8=\begin{pmatrix}1&1\\0&1
\end{pmatrix}\in N$ and hence $\OO=\G_\8\sm N\subset X$.  We denote by $dn$ the $N$-invariant measure on $\OO$ coming from $N$ such that $vol_{dn}(\OO)=1$. 

Let $\psi_n(n_x)=e^{2\pi inx}$ be a {\it non-trivial} character $\psi_n:N\to\CC^\times$ trivial on $\G_\8$. We can view such a character as a function on $\OO$ and $\{\psi_n\}_{n\in\ZZ}$ as an orthonormal system in $L^2(\OO,dn)$. Consider the Whittaker  functional defining the $\psi_n$'th Fourier coefficient on the automorphic representation $(\pi,\nu)$
$$\WW^{\rm aut}_n:\,V_{\pi}\rightarrow\mathbb{C},
\quad v\mapsto\int_{\OO}\phi_v(o)\bar\psi_n(o)do=\int_0^1
\phi_v(n_x)\bar\psi_n(n_x)dx\ .$$
The functional $\WW^{\rm aut}_n\in{\rm Hom}_N(V_{\pi},{\psi}_n)$ is  ${\psi}_n$-equivariant, and  gives the $n$-th Fourier coefficient $a_n(v)$ in the expansion of $\phi_v$ restricted to $\OO$: $\phi_v|_{_{\OO}}=\sum_na_n(v){\psi}_n$. 
It is known that the space ${\rm Hom}_N(V_{\pi},{\psi}_n)$  is at most one dimensional if $\psi_n$ is non-trivial. We will assume that $\pi$ is a representation of the principal series and moreover that $(\pi,\nu)$ is {\it cuspidal}, that is, $\WW^{\rm aut}_0\equiv 0$.  We define the  model Whittaker functional in the linear model $V_{\lambda}$ of $\pi$ by:
$$\WW^{\rm mod}_n:\,V_{\lambda}\rightarrow\mathbb{C},\quad v\mapsto\langle v,\psi_n\rangle=\tfrac{1}{\pi}\int_{\mathbb{R}}v(x)
\overline{\psi}_n(x)dx\ ,$$ for any $n\not=0$. The last integral is not absolutely convergent, but easily could be extended from compactly supported functions to smooth vectors in $V_\lm$ (i.e., decaying as $|x|^{\lm-1}$ at infinity).
We have $\WW^{\rm mod}\in{\rm Hom}_N(V_{\pi},{\psi}_n)$ and hence  there exists a scalar $c_n=c_n(\nu)\in\mathbb{C}$ such that $$\WW^{\rm aut}_n(\phi_v)=c_n\cdot
\WW^{\rm mod}_n(v)$$ 
for any $v\in V_{\lambda}$, $n\ne0$. The Plancherel formula in $L^2(\OO,dn)$  for the vector $\pi(g)v$ gives 
\begin{equation}\label{eq-2}
||\phi_{\pi(g)v}|_\OO||^2_{L^2(\OO)}=\sum_n|c_n|^2\big|\WW^{\rm mod}_n\big(\pi(g)v\big)\big|^2\ .
\end{equation}

We now want to integrate this identity with respect to $g\in G$. Clearly action of $N$ on the left is trivial and we can  consider $g\in N\sm G\simeq AK$. Moreover, we will apply this to a $K$-type vector $v=e_m\in V_\pi$ and hence can consider $g\in A$.  
Let $\al\in C^\8_c(N\sm G/K)$ be a smooth  function on $G$ which is left $N$-invariant, right $K$-invariant and is compactly supported on the quotient $N\sm G/K$. The exact choice of $\al$ will be made  later. Let $\tau([a])dnd[a]dk$ be the invariant measure on $G$ written in the Iwasawa $NAK$ coordinates with $\tau(a)=\tau([a])=a^{-2}$  and $d[a]=da/a$ for $[a]={\rm diag}(a,a\inv)$.  
We obtain for a $K$-type $v$, 
\begin{eqnarray*}
\int_A\int_\OO\big|\phi_{\pi(a)v}\big)\big|^2\al(a)\tau(a)d[a]
&=&
\int_K\int_A\left(\int_{\OO}\big|
\phi_v(oak)\big|^2dx\right)\al(ak)\tau(a)dod[a]dk\\[0.2cm]
&\leq&\int_K\int_A
\int_{\G_\8\sm N}\big|
\phi_v(nak)\big|^2\al(nak)\tau(a)dnd[a]dk\\[0.2cm]
&\overset{(\#)}{=}&\int_{\G_\8\sm G}\big|
\phi_v(g)\big|^2\al(g)
dg
=\int_{X}
\big|
\phi_v(x)\big|^2H_\al(x)dx\ ,
\end{eqnarray*}
where $H_\al(g)=\sum_{\g\in\G_\8\sm \G}\al(\g g)$ (note that we assume that $\OO=\G_\8\sm N=\G N\subset X$).   At the step (\#), we have used the integral formula for functions on  $G$  written in the $NAK$ decomposition (see \cite{kn}). On the right in (\ref{eq-2}), we obtain  
$$
\sum_n|c_n|^2\int_A\big|\langle\pi_{\lambda}(a)v,\WW^{\rm mod}_n\rangle\big|^2\al(a)\tau(a)d[a].$$
Thus we have
\begin{equation*}
\int_{X}
\big|
\phi_v(x)\big|^2H_\al(x)dx=\sum_n|c_n|^2\int_A\big|\langle\pi_{\lambda}(a)v,\WW^{\rm mod}_n\rangle\big|^2\al(a)a^{-2}\frac{da}{a}\ 
\end{equation*}
for any $v\in V_{\lambda}$.

\subsection{Airy asymptotic}
For $[a]={\rm diag}(a,a\inv)\in A$ (with $a\leq 1$), we have 
\begin{eqnarray}\label{Airy-uni}&& \langle\pi_{\lambda}([a])e_m,\WW^{\rm mod}_n\rangle=\\ &&\langle e_m,\pi_{\lambda}\big([a]^{-1}\big)\WW^{\rm mod}_n\rangle
={\pi}^{-\haf}|a|^{1-\lm}\int_{\RR}
e^{im\arctan x-in(a^2x)}(1+x^2)^{\tfrac{\lambda-1}{2}}dx\ .\nonumber
\end{eqnarray}
The integral \eqref{Airy-uni} has  Airy type critical point at $x=0$. We choose the partition of unity $\chi_0$, $\chi_\8$ as before and study the contribution near $x=0$ (it is easy to see that the contribution away from the interval $[-1,1]$ is negligible). We introduce parameters $1>\bt$, $1>\eps>0$ satisfying $a=\beta(1+\eps)$ and $m=-\bt^2n$. Let $S(\bt,\eps;x)=\arctan x+(1+\eps)^2x$ be the phase function and $f(x)=(1+x^2)^{\tfrac{\lambda-1}{2}}$ be the amplitude  in the integral in \eqref{Airy-uni}
\begin{equation}\label{uni-int}J_{m,n}(a)=\int_{\RR}
e^{imS(\bt,\eps;x)}f(x)\chi_0(x)dx\ .
\end{equation}   We have then  $S(\bt,\eps;x)=2\eps x-x^3/3+({\rm  smaller \ order\ terms \ in \ } x, \eps, \bt)$. This implies  that the value of the integral is well-approximated by $(n\bt^2)^{-1/3}{\rm Ai}((n\bt^2)^{2/3}\eps)$, and hence  the integral in \eqref{Airy-uni} is bounded from below  by $c\cdot (n\bt^2)^{-1/3}$ for some explicit constant $c>0$, as long as the following Airy condition holds:
\begin{equation}\label{cond-uni}
0\leq \eps \leq (n\bt^2)^{-2/3}\leq 1/2. 
\end{equation} Finally, we get the lower bound for the matrix coefficient of the form 
\begin{equation}\label{mat-c-bound-uni1}
|\langle\pi_{\lambda}(a)e_m,
\WW^{\rm mod}_n\rangle|^2 =\pi\inv|a|^{2}|J_{m,n}(a)|^2\geq C\cdot n^{-2/3}\beta^{2/3}\ ,
\end{equation} for some explicit constant $C>0$ depending on $\lm$, but not on $a$, $m$ and $n$.

We obtain the following
\begin{proposition}\label{m-c-bound-prop-uni} Let $a\leq 1$, $m, n>0$, be such that parameters $\beta=(m/n)^\haf\leq 1$ and $\eps=a/\beta -1$ satisfy the Airy condition \eqref{cond-uni}.  We have then the following lower bound
\begin{equation}\label{mat-c-bound-uni}
|\langle\pi_{\lambda}(a)e_m,
\WW^{\rm mod}_n\rangle|^2 \geq C\cdot n^{-2/3}\beta^{2/3}\ ,
\end{equation} for some explicit constant $C>0$ depending on $\lm$, but not on $a$ and $n$.
\end{proposition}

\subsubsection{Lattice points counting for unipotent periods} We briefly discuss the lattice points counting problem in order to show that an effective version of it gives rise to a non-trivial bound on the Fourier coefficient of a fixed Maass cusp form. We do not try to optimize the exponent in the bound since we do not see how our method could compete with best available bounds.

We have to analyze the behavior of the function $H_{\al_{\bt,\eps}}$ for $\bt$ and $\eps$ satisfying conditions of Proposition \ref{m-c-bound-prop-uni}. As in previous sections, we can reduce the question of bounding $H_{\al_{\bt,\eps}}$ to the question of obtaining effective lattice points counting for $\G_\8\sm \G$ in the following set
$B_T=\G_\8\sm NA_TK$, where $A_T=\{ [a]=diag(a,a\inv)\ | \ 1\geq  a\geq T\inv\}$. This a well-known problem of counting lattice points in sectors (see  \cite{gos}). 
The resulting bound is non-uniform on $X$ (i.e., the function $H_{\al_{\bt,\eps}}$ growth in the cusp as $\bt\to 0$). We have to apply the ``cutting'' procedure from Section  \ref{non-comp} in order to modify the function $H_{\al_{\bt,\eps}}$. 

It is also possible to apply the spectral  approach of Selberg in order to bound the function $H_{\al_{\bt,\eps}}$. Since the function $\al_{\bt,\eps}$ (or the characteristic function of $B_T$) is $N$-invariant,  its spectral decomposition  involves only the Eisenstein series. From this, one can deduce a bound on $||H_{\al_{\bt,\eps}}||_{L^2(X)}$ and then apply the bound on $L^4$-norm of $K$-types from \cite{br}. In particular, the best possible bound on the average size of the number of lattice points in sectors which is on par with the bound for the spherical counting obtained in \cite{hp} would give $H_{\al_{\bt,\eps}}(x)=\int \al_{\bt,\eps}(a)\tau(a)d[a]\cdot {\bf 1}+R_{\bt,\eps}(x)$ with $||R_{\al_{\bt,\eps}}||_{L^2(X)}\leq c \bt\inv$ and $\int \al_{\bt,\eps}(a)\tau(a)d[a]\sim \bt^{-2}\eps$.

This in turn implies the following

\begin{theorem} There exists a constant $C>0$  depending on $\Gamma$ and $\nu$  such that 
$|c_n|^2\leq C{\cdot}|n|^{6/7}$.
\end{theorem}

\vskip 1cm
\end{document}